\documentclass[a4 paper, 12pt]{article} 
\usepackage{amssymb,amsthm}
\usepackage{amsmath}

\usepackage{tikz}
\usetikzlibrary{decorations.pathmorphing}
\tikzset{snake it/.style={decorate, decoration=snake}}

\newcommand{\projecttitle}{X-ray Compton scattering tomography}
\newcommand{\projectauthor}{By James Webber}
\theoremstyle{plain}
\newtheorem{theorem}{Theorem}
\newtheorem{lemma}{Lemma}
\theoremstyle{definition}

\newtheorem{corollary}{Corollary}

\def\titlep{\thispagestyle{empty}\null\vskip1in\begin{center}
           \Huge\uppercase{\projecttitle}
\end{center}
\vspace{10mm}
\begin{center}
\text{\projectauthor}
\end{center}
\vspace{10mm}
}

\begin{document}
\titlep
\begin{abstract}
We lay the foundations for a new fast method to reconstruct the electron density in x-ray scanning applications using measurements in the dark field. This approach is applied to a type of machine configuration with fixed energy sensitive (or resolving) detectors, and where the X-ray source is polychromatic.
We consider the case where the measurements in the dark field are dominated by the Compton scattering process. This leads us to a 2D inverse problem where we aim to reconstruct an electron density slice from its integrals over discs whose boundaries intersect the given source point. We show that a unique solution exists for smooth densities compactly supported on an annulus centred at the source point.

Using Sobolev space estimates we determine a measure for the ill posedness of our problem based on the criterion given by Natterer in \cite{nat1}. In addition, 
with a combination of our method and the more common attenuation coefficient reconstruction, we show under certain assumptions that the atomic number of the target is uniquely determined. 

We test our method on simulated data sets with varying levels  of added pseudo random noise.
\end{abstract}

\section{Introduction}
In this paper we investigate the potential for the use of incoherent scattered data for 2D reconstruction in x-ray scanning applications. The use of scattered data for image reconstruction is considered in the literature, typically for applications in gamma ray imaging, where the photon source is monochromatic \cite{pal1, cone, norton}. However, in many applications (e.g security screening of baggage) a type of x-ray tube is often used that generates a polychromatic spectrum of initial photon energies (see section \ref{phys} for an example spectrum). There has been recent interest in the use of energy sensitive detectors in tomography \cite{pca1, pca2}, and in the present paper their application is key to the ideas presented.

Our main goal is to show that the electron density may be reconstructed analytically using the incoherent scattered data and to lay the foundations for a practical reconstruction method based on our theory. We apply our method to a machine configuration commonly used in x-ray CT.
In addition, by use of the reconstructed density values in conjunction with an attenuation coefficient reconstruction, we show under the right assumptions that the atomic number of the target is uniquely determined. 

For a photon incident upon an electron Compton (incoherently) scattering at an angle $\omega$ with initial energy $E_{\lambda}$, the scattered energy $E_s$ is given by the equation:
\begin{equation}
\label{equ1}
E_s=\frac{E_{\lambda}}{1+\left(E_{\lambda}/E_0\right)\left(1-\cos\omega\right)}
\end{equation}
where $E_0\approx 511$keV is the electron rest energy. Equation (\ref{equ1}) implies that $\omega$ remains fixed for any given $E_s$ and $E_{\lambda}$. So in the case of a monochromatic source, assuming only single scatter events, for every fixed measured energy $E_s$ (possible to measure if the detectors are energy-resolved) the locus of scattering points is a circular arc intersecting the source and detector in question. For example, refer to \cite{pal1,norton}.

In an x-ray tube a cathode is negatively charged and electrons are accelerated by a large voltage ($E_{\text{max}}$ kV) towards a positively charged target material (e.g Tungsten). A small proportion of the initial electron energy ($\approx 1\%$) is converted to produce photons. Due to energy conservation, the resulting photon energies are no more than $E_{\text{max}}$ keV. 
So in the polychromatic source case, again assuming only single scatter events, for each given data set (photon intensity recorded with energy $E_s$), the set of scatterers lie on a collection of circular arcs intersecting the source and detector points. Together these form a toric section in which the photons scatter, with a maximum scattering angle $\omega_{\text{max}}$ given by:
\begin{equation}
\cos\left(\omega_{\text{max}}\right)=1-\frac{E_0\left(E_{\text{max}}-E_s\right)}{E_{s}E_{\text{max}}}
\end{equation}
See figure \ref{fig1} below:
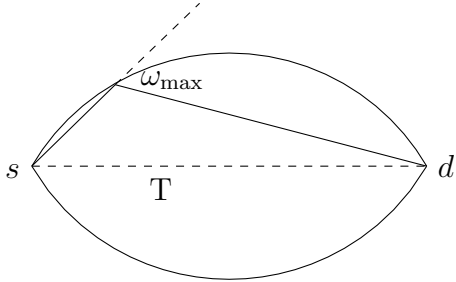
\begin{figure}[!h]
\centering
\begin{tikzpicture}[scale=3]
\path ({-cos(30)},0.5) coordinate (S);
\path ({cos(30)}, 0.5) coordinate (D);
\path (-0.5,0.86) coordinate (w);
\path (-0.13,1.22) coordinate (a);
\draw [domain=30:150] plot ({cos(\x)}, {sin(\x)});
\draw [thin, dashed] (S) -- (D);
\draw (S) -- (w);
\draw (w) -- (D);
\draw [thin,dashed] (w)--(a);
\node at (-0.95,0.48) {$s$};
\node at (0.95,0.5) {$d$};
\node at (-0.25,0.88) {$\omega_{\text{max}}$};
\node at (-0.3,0.4) {T};
\draw [domain=-0.866:0.866] plot(\x, {-sqrt(1-pow(\x,2))+1});
\end{tikzpicture}
\caption{A toric section T  in which the photons scatter with tips at source and detector points $s$ and $d$.}
\label{fig1}
\end{figure}
\\In the present paper we consider a setup consisting of a ring of fixed energy sensitive detectors and a single rotating fan beam polychromatic source. See figure \ref{fig3}. With this setup we can measure photon intensity in the dark field. We image an electron density $f:\mathbb{R}^2 \to \mathbb{R}$ compactly supported within the detector ring (the blue and green circle in figure \ref{fig3}), with $f\geq 0$.
\begin{figure}[!h]
\centering
\begin{tikzpicture}[scale=3.25]
\draw [ blue, dashed] (0,0) circle [radius=0.7];;
\draw [red, very thin] (-0.3,-0.2) rectangle (0.3,0.2);
\draw [very thin] (0.651,-0.546) circle [radius=0.8];
\node at (-0.05,-1.08) {$s$};
\node at (0.7,0.34) {$d$};
\node at (0.651+0.8,0) {$D$};
\node at (-0.15,0.1) {$f$};
\node at (0.2,-0.4) {R};
\draw [very thin, dashed] (0,-1)--(0,0.7);
\draw [very thin, dashed] (0,-1)--(0.1,0.6928);
\draw [very thin, dashed] (0,-1)--(-0.1,0.6928);
\draw [very thin, dashed] (0,-1)--(0.2,0.6708);
\draw [very thin, dashed] (0,-1)--(-0.2,0.6708);
\draw [very thin, dashed] (0,-1)--(0.3,0.6324);
\draw [very thin, dashed] (0,-1)--(-0.3,0.6324);
\draw [very thin, dashed] (0,-1)--(0.4,0.5744);
\draw [very thin, dashed] (0,-1)--(-0.4,0.5744);
\draw [very thin] (0.9,0.9)--(0.3,0.6324);
\node at (1.1,1) {light field};
\draw [very thin] (-0.9,-0.9)--(-0.3,-0.6324);
\node at (-1.1,-1) {dark field};
\draw[snake it, very thin] (0.03,-0.1)->(0.65,0.255);
\draw[snake it, very thin] (-0.05,-0.4)->(0.65,0.255);
\draw [thick, green, domain=-0.4:0.4] plot(\x, {sqrt(0.49-pow(\x,2))});
\draw [<-,domain=-0.45:0] plot(\x, {-sqrt(1-pow(\x,2))});
\end{tikzpicture}
\caption{An example machine configuration is displayed. A disc $D$ whose boundary intersects a point source $s$ and a detector $d$ forms the scattering region R. The source $s$ travels along the circular path shown. Detectors under direct exposure to the initial x-ray beam are said to be in the light field. Detectors not in the light field are said to be in the dark field.}
\label{fig3}
\end{figure}
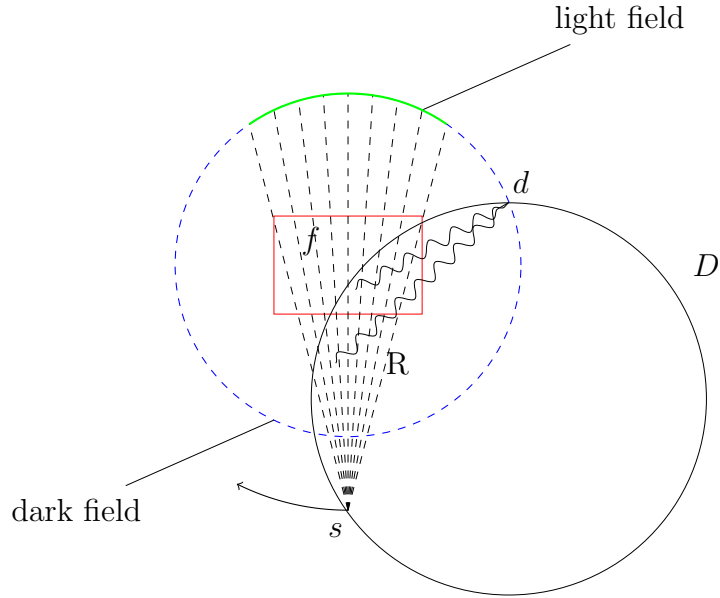
If we assume an equal scattering probability throughout the region $\text{R}=D\cap \text{supp}\left(f\right)$ leaving only the electron density to vary, and if we assume that the majority of scattering events occur within R, then in this case the integral of $f$ over $D$ is approximately determined by the scattered intensity recorded at the detector $d$ with some fixed energy $E_s$. See the appendix for an example application where these approximations are valid. With these assumptions and with suitable restrictions on the support of $f$, we aim to reconstruct $f$ from its integrals over discs whose boundaries intersect a fixed point, namely the source at a given position along its scanning path. 

In section \ref{sec1}, we present a disc transform and go on to prove our main theorem (Theorem \ref{th1}), which explains the relationship between our transform and the straight line Radon transform. As a corollary to this theorem, with known results on the Radon transform, we show that a unique solution exists on the domain of smooth functions compactly supported on an annulus centred at the origin. Here based on the criterion of Natterer in \cite{nat1} and using the theory of Sobolev space estimates, we determine a measure for the ill posedness of our problem.

In section \ref{phys}, we discuss a possible means to approximate the physical processes such as to allow for the proposed reconstruction method. Here we also present a least squares fit for the total cross section (scattering plus absorbtion) in terms of $Z$ (the atomic number). From this, we show that $Z$ is uniquely determined by the attenuation coefficient and electron density.

In section \ref{res} we apply our reconstruction formulae to simulated data sets, with varying levels of added pseudo random noise. This is applied to the given machine configuration. We recover a simple water bottle cross section image (a circular region of uniform density 1) and reconstruct the atomic number in each case using the curve fit presented in section \ref{phys}. To give an example reconstruction of a target not of uniform density, we also present reconstructions of a simulated hollow tube cross section.

\section{A disc transform}
\label{sec1}
In this section we aim to recover a smooth function compactly supported on an annulus centred at the origin $O$ from its integrals over discs whose boundaries intersect $O$ (the given source position).

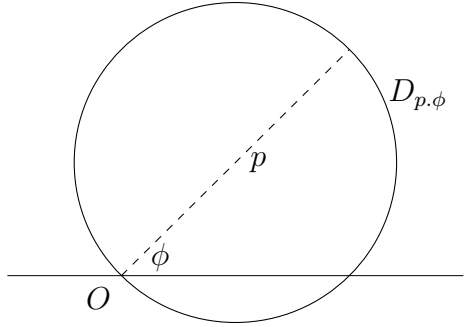
\begin{figure}[!h]
\centering
\begin{tikzpicture}[scale=3]
\draw [very thin] (-0.5,0)--(1.5,0);
\draw (0.5,0.5) circle [radius=0.707];
\draw [thin, dashed] (0,0)--(1,1);
\node at (0.6,0.5) {$p$};
\node at (0.17,0.08) {$\phi$};
\node at (-0.1,-0.1) {$O$};
\node at (1.3,0.8) {$D_{p.\phi}$};
\end{tikzpicture}
\label{figure8}
\caption{A disc $D_{p,\phi}$ with its boundary intersecting $O$.}
\end{figure}

Let $D_{p,\phi}$ denote the set of points on the disc whose boundary intersects the origin, with centre given in polar coordinates as $\left(p/2,\phi\right)$. See figure \ref{figure8}. Let $C^{\infty}\left(\Omega\right)$ be the set of smooth functions on $\Omega \subseteq \mathbb{R}^n$ and let $C_0 ^{\infty}\left(\Omega\right)$ denote the set of smooth functions compactly supported on $\Omega$. Let $Z^{+}=\mathbb{R}^{+}\times S^1$ and for a function in the plane $f:\mathbb{R}^2\to \mathbb{R}$, let $F:Z^{+} \to \mathbb{R}$ be defined as $F\left(\rho,\theta\right)=f\left(\rho \cos \theta,\rho \sin \theta\right)$. Then we define the disc transform $\mathcal{D}_1 : C_{0} ^{\infty}\left(\mathbb{R}^2\right)\to C^{\infty}\left(Z^{+}\right)$ as:
\begin{equation}
\begin{split}
\label{equdef1}
\mathcal{D}_1f\left(p,\phi\right)&=\iint_{D_{\frac{1}{p},\phi}}f \mathrm{d}A=\int_{-\frac{\pi}{2}}^{\frac{\pi}{2}}\int_{0}^{\frac{\cos{\theta}}{p}}\rho F\left(\rho,\theta+\phi\right)\mathrm{d}\rho \mathrm{d}\theta
\end{split}
\end{equation}

After making the change of variables:
\begin{equation}\rho=r\cos\psi, \ \ \ \theta=\psi,\ \ \ \mathrm{d}\rho \mathrm{d}\theta=\cos\psi \mathrm{d}r \mathrm{d}\psi\end{equation}
in equation (\ref{equdef1}), we have:
\begin{equation}\mathcal{D}_1 f\left(p,\phi\right)=\int_{0}^{\frac{1}{p}} \int_{-\frac{\pi}{2}}^{\frac{\pi}{2}}r\cos^{2}\psi F\left(r\cos \psi,\psi+\phi\right)\mathrm{d}\psi \mathrm{d}r\end{equation}
We now present further definitions which will be important in the following subsection (section \ref{sob}), where we provide our Sobolev space estimates. Let $Z=\mathbb{R}\times S^1$ denote the unit cylinder in $\mathbb{R}^3$. Then we define $\mathcal{D}_2 : C_{0}^{\infty}\left(\mathbb{R}^2\right) \to C^{\infty}\left(Z\right)$ as follows:
\begin{equation}\mathcal{D}_2 f\left(p,\phi\right)=\left\{
	\begin{array}{ll}
		\mathcal{D}_1f\left(p,\phi\right) & \mbox{$p>0$}\\
                    \frac{\mathcal{D}_1f\left(0^{+},\phi\right)+\mathcal{D}_1f\left(0^{+},\phi+\pi\right)}{2} &  \mbox{$p=0$}\\
                     \mathcal{D}_1f\left(-p,\phi+\pi\right) &  \mbox{$p<0$}
           \end{array}
\right.\end{equation}
which is piecewise continuous as a function of $p$. We can remove this discontinuity by adding the function:
\begin{equation}
c\left(\phi\right)\text{sgn}\left(p\right)=\left\{
	\begin{array}{ll}
		c\left(\phi\right) & \mbox{$p>0$}\\
                    0 &  \mbox{$p=0$}\\
                     -c\left(\phi\right) &  \mbox{$p<0$}
           \end{array}
\right.
\end{equation}
where $c\left(\phi\right)=  \frac{\mathcal{D}_1f\left(0^{+},\phi+\pi\right)-\mathcal{D}_1f\left(0^{+},\phi\right)}{2}$. We define $\mathcal{D} : C_{0}^{\infty}\left(\mathbb{R}^2\right) \to C^{\infty}\left(Z\right)$ as:
\begin{equation}\mathcal{D}f\left(p,\phi\right)=\mathcal{D}_2f\left(p,\phi\right)+c\left(\phi\right)\text{sgn}\left(p\right)\end{equation}
Let $L_{p.\phi}=\{\left(x,y\right)\in \mathbb{R}^2 : x\cos \phi +y\sin \phi =p\}$ be the set of points on a line. Then we define the Radon transform $R: C_{0}^{\infty}\left(\mathbb{R}^2\right) \to C^{\infty}\left(Z\right)$ as:
\begin{equation}Rf\left(p,\phi\right)=\int_{L_{p.\phi}}f\mathrm{d}s\end{equation}
We are now in a position to prove our main theorem, where we give the explicit relation between $\mathcal{D}$ and the Radon transform $R$ for smooth functions on an annulus.

\begin{theorem}
\label{th1}
Let $A_{r_1,r_2}=\{x\in \mathbb{R}^2 : r_1<|x|<r_2\}$ be the annulus centred on $O$ with inner radius $r_1>0$ and outer radius $r_2$. Let $f\in C_0^{\infty}\left(A_{1,r}\right)$ for some $r>1$ and let $\tilde{f}\in C_0^{\infty}\left(A_{1/r,1}\right)$ be defined as $\tilde{f}\left(x\right)=\frac{1}{|x|^4}f\left(\frac{x}{|x|^2}\right)$. Then $\frac{\partial}{\partial p}\mathcal{D}f=-R\tilde{f}$.
\begin{proof}
Let $\tilde{F}$ and $F$ be defined as $\tilde{F}\left(\rho,\theta\right)=\tilde{f}\left(\rho \cos\theta,\rho\sin\theta\right)$ and ${F}\left(\rho,\theta\right)={f}\left(\rho \cos\theta,\rho\sin\theta\right)$. Then from our definition of $\tilde{f}$, we have $\tilde{F}\left(\rho,\theta\right)=\frac{1}{\rho^4}F\left(\frac{1}{\rho},\theta\right)$. Now we have:
\begin{equation}
\begin{split}
\frac{\partial}{\partial p}\mathcal{D}_1f\left(p,\phi\right)&=-\frac{1}{p^3} \int_{-\frac{\pi}{2}}^{\frac{\pi}{2}}\cos^{2}\psi F\left(\frac{\cos \psi}{p},\psi+\phi\right)\mathrm{d}\psi \\
&= -p\int_{-\frac{\pi}{2}}^{\frac{\pi}{2}}\tilde{F}\left(\frac{p}{\cos \psi},\psi+\phi\right)\frac{\mathrm{d}\psi }{\cos^2 \psi}\\
&=-R\tilde{f} \left(p,\phi\right)\ \ \ \text{for}\ \ \ p\geq 0
\end{split}
\end{equation}
and hence $\frac{\partial}{\partial p}\mathcal{D}f\left(0,\phi\right)=\frac{\partial}{\partial p}\mathcal{D}_1f\left(0,\phi\right)=\frac{\partial}{\partial p}\mathcal{D}_1f\left(0,\phi+\pi\right)$. So the partial derivative of $\mathcal{D}f$ with respect to $p$ exists and is continuous for all $p\in \mathbb{R}$, and $\frac{\partial}{\partial p}\mathcal{D}f=-R\tilde{f}$.
\end{proof}
\end{theorem}
We now aim to prove injectivity of the disc transform $\mathcal{D}$ on the domain of smooth functions compactly supported on an annulus. First we state Helgason's support theorem \cite{he}.

\begin{theorem}
Let $X$ be a compact convex set in $\mathbb{R}^n$ and let $f$ be continuous on $\mathbb{R}^n/X$. If $Rf=0$ for all $p$ and $\phi$ such that $L_{p,\phi} \cap X=\emptyset$ and $f$ is rapidly decreasing, in the sense that:
\begin{equation}|x|^k f\left(x\right) \to 0 \ \ \ \text{as} \ \ \ |x| \to \infty \ \ \ \forall k\in \mathbb{N}\end{equation}
then $f\left(x\right)=0$ for all $x\notin X$. 
\end{theorem}

\begin{corollary}
\label{cor1}
Let $f\in C_0 ^{\infty}\left(A_{1,r}\right)$ for some $r>1$, and let $Z_r=\{(p,\phi)\in Z : 1/r< p< 1, \phi \in [0,2\pi]\}$. Then $f$ is uniquely determined by $\mathcal{D}f$ known for all $(p,\phi) \in Z_r$.

\begin{proof}
Let:
\begin{equation}
f\in \{ f\in  C_0 ^{\infty}\left(A_{1,r}\right): \mathcal{D}f=0\ \ \text{for all}\ \ (p,\phi)\in Z_r\}
\end{equation}
and let $\tilde{f}$ be defined as in Theorem \ref{th1}. Then by Theorem \ref{th1}, we have:
\begin{equation}
\label{unique}
\tilde{f}\in \{f\in C_0^{\infty}\left(A_{1/r,1}\right) : Rf=0 \ \ \text{for all}\ \ (p,\phi)\in Z_r\}
\end{equation}
and hence $\tilde{f}$ is rapidly decreasing. Let $X=\{x\in\mathbb{R}^2 : |x|\leq1/r\}$. Then $X$ is clearly compact and convex. By (\ref{unique}), $R\tilde{f}=0$ for all $p$ and $\phi$ such that $L_{p,\phi}\cap X=\emptyset$. So by Helgason's support theorem, we have that $\tilde{f}(x)=0$ for all $x\notin X$. The result follows. 
\end{proof}
\end{corollary}
For the proposed machine configuration, we can define the set of points within the detector ring formally as $D_r=\{\left(x,y\right)\in \mathbb{R}^2 : x^2+\left(y-\left(r+1\right)/2\right)^2<\left(r-1\right)^2/4\}$, where $r>1$ depends on the machine specifications (i.e the detector ring radius and the source path radius). We now have:
\begin{corollary}
\label{cor2}
Let $f\in C_0^{\infty}\left(D_r\right)$. Let $\partial D_r$ denote the boundary of $D_r$, and let $\text{R}_{p,\phi}=D_{\frac{1}{p},\phi}\cap D_r$. Then the values of $\mathcal{D}_1f$ for $p$ and $\phi$ such that:
\begin{equation}
\label{equcor2}
\text{R}_{p,\phi}\neq \emptyset\ \  \text{and}\ \ \partial D_{\frac{1}{p},\phi}\cap \partial D_r\neq \emptyset
\end{equation}
 determine $f$ uniquely.
\begin{proof}
We consider two cases. If $D_r\subset D_{\frac{1}{p},\phi}$ then $\mathcal{D}_1f\left(p,\phi\right)=\mathcal{D}_1f\left(\frac{1}{r},\pi/2\right)$, which is known as condition (\ref{equcor2}) is satisfied for $p=1/r$ and $\phi=\pi/2$. If $D_r\cap D_{\frac{1}{p},\phi}=\emptyset$, then $\mathcal{D}_1f\left(p,\phi\right)=0$. In any other case, $\mathcal{D}_1f$ is known by our assumption. Hence we have a full data set for $\mathcal{D}_1f$ and hence for $\mathcal{D}f$. The result follows from Corollary \ref{cor1}.
\end{proof}
\end{corollary}
So for the proposed application, we see from the above corollaries that for any given source position, the incoherent scattered data is sufficient to reconstruct the target density uniquely. 

\subsection{Sobolev space estimates}
\label{sob}
In this section we provide Sobolev space estimates for the disc transform $\mathcal{D}$. From these we obtain an upper bound for the least squares error in our solution $f$ in terms of $\epsilon$, where $\epsilon$ is an upper bound for the least squares error in our measurements. First we define our Sobolev spaces and the norms which will be used in our estimates.

Let $\Omega \subset \mathbb{R}^n$ be an arbitrary domain and let $L^2\left(\Omega\right)$ denote the set of square integrable functions on $\Omega$. 
We define the Fourier transform of a function $f\in L^2\left(\mathbb{R}^n\right)$ as:
\begin{equation}\hat{f}\left(\xi\right)=\left(2\pi\right)^{-n/2}\int_{\mathbb{R}^n}f\left(x\right)e^{-ix\cdot \xi}\mathrm{d}x\end{equation}
Then we can define Sobolev spaces $H^{\alpha}\left(\mathbb{R}^n\right)$ of real degree $\alpha\in \mathbb{R}$ as:
\begin{equation}H^{\alpha}\left(\mathbb{R}^n\right)=\{\text{tempered distributions}\ f : \left(1+|\xi|^2\right)^{\alpha/2}\hat{f}\left(\xi\right)\in L^2\left(\mathbb{R}^n\right)\}\end{equation}
with the norm:
\begin{equation}
\label{norm2}
{\|f\|^2_{{H}^{\alpha}\left(\mathbb{R}^n\right)}}=\int_{\mathbb{R}^n}\left(1+|\xi|^2\right)^{\alpha}|\hat{f}\left(\xi\right)|^2\mathrm{d}\xi
\end{equation}
For functions on the cylinder $Z\subset \mathbb{R}^3$, we have the norm:
\begin{equation}
{\|f\|^2_{{H}^{\alpha}\left(Z\right)}}=\int_{S^1}\int_{\mathbb{R}}\left(1+\sigma^2\right)^{\alpha}|\hat{f}\left(\sigma,\phi\right)|^2\mathrm{d}\sigma\mathrm{d}\phi
\end{equation}
where the Fourier transform of $Rf$ is taken with respect to the variable $p\in \mathbb{R}$. We now state some preliminary results on Sobolev spaces and the Radon transform which will be used in our theorems \cite[pages 11 and 203]{nat1}.

\begin{theorem}
For $f\in S\left(\mathbb{R}^2\right)$, where $S\left(\mathbb{R}^2\right)$ is the Schwartz space on $\mathbb{R}^2$, we have:
\begin{equation}\hat{Rf}\left(\sigma,\phi\right)=\left(2\pi\right)^{\left(-1/2\right)}\hat{f}\left(\sigma\Phi\right),\ \ \sigma\in \mathbb{R}\end{equation}
where $\Phi=\left(\cos \phi,\sin\phi\right)$ and the Fourier transform of $Rf$ is taken with respect to the $p$ variable only.
\end{theorem}

\begin{theorem}
Let $k=\left(k_1,\ldots,k_n\right)$ be some multi index and let $D^{k}=\frac{\partial^{k_1}}{\partial x_1^{k_1}}\cdots\frac{\partial^{k_n}}{\partial x_n^{k_n}}$, where the $\frac{\partial}{\partial x_i}$ are defined in the weak sense. Let $m$ be a positive integer and let $\sigma \in \left(0,1\right)$. Then for $\alpha=m+\sigma$, the norm (\ref{norm2}) is equivalent to the norm:
\begin{equation}
{\|f\|^2_{{H}^{\alpha}\left(\Omega\right)}}={\|f\|^2_{H^{m}\left(\Omega\right)}}+\sum_{|k|=m}\iint_{\Omega\times \Omega}\frac{|D^kf\left(x\right)-D^kf\left(y\right)|^2}{|x-y|^{n+2\sigma}}\mathrm{d}x\mathrm{d}y
\end{equation}
when $\Omega=\mathbb{R}^n$.
\end{theorem}

We now prove a slice theorem for the disc transform $\mathcal{D}$.

\begin{lemma}
Let $f\in C^{\infty}_0\left(A_{1,r}\right)$ for some $r>1$ and let $\tilde{f}$ be defined as in Theorem \ref{th1}. Then we have:
\begin{equation}-i\sigma\hat{\mathcal{D}f}\left(\sigma,\phi\right)=\left(2\pi\right)^{\left(-1/2\right)}\hat{\tilde{f}}\left(\sigma\Phi\right),\ \ \sigma\in \mathbb{R}\end{equation}
where $\Phi=\left(\cos \phi,\sin\phi\right)$ and the Fourier transform of $\mathcal{D}f$ is taken with respect to the $p$ variable.

\begin{proof}
 Let $\mathcal{D}_2$ and $c\left(\phi\right)$ be as defined in section \ref{sec1}. Then we have:
\begin{equation}
\begin{split}
\hat{R\tilde{f}}\left(\sigma,\phi\right)&=\left(2\pi\right)^{-1/2}\int_{-\infty}^{\infty}R\tilde{f}\left(p,\phi\right)e^{-ip\sigma}\mathrm{d}p\\
&=\left(2\pi\right)^{-1/2}\Big[\int_{-\infty}^{0}R\tilde{f}\left(p,\phi\right)e^{-ip\sigma}\mathrm{d}p+\int_{0}^{\infty}R\tilde{f}\left(p,\phi\right)e^{-ip\sigma}\mathrm{d}p\Big]\\
&=-\left(2\pi\right)^{-1/2}\Big[\int_{-\infty}^{0}\frac{\partial}{\partial p}\mathcal{D}_2f\left(p,\phi\right)e^{-ip\sigma}\mathrm{d}p+\int_{0}^{\infty}\frac{\partial}{\partial p}\mathcal{D}_2f\left(p,\phi\right)e^{-ip\sigma}\mathrm{d}p\Big]\\
&=-\left(2\pi\right)^{-1/2}\left(\mathcal{D}_2f\left(0^{-},\phi\right)-\mathcal{D}_2f\left(0^{+},\phi\right)\right)-\frac{i\sigma}{\left(2\pi\right)^{1/2}}\int_{-\infty}^{\infty}\mathcal{D}_2f\left(p,\phi\right)e^{-ip\sigma}\mathrm{d}p\\
&=-\frac{2c\left(\phi\right)}{\left(2\pi\right)^{1/2}}-i\sigma\hat{\mathcal{D}_2f}\left(\sigma,\phi\right)\\
&=-i\sigma\big[\frac{2c\left(\phi\right)}{i\sigma\left(2\pi\right)^{1/2}}+\hat{\mathcal{D}_2f}\left(\sigma,\phi\right)\big]\\
&=-i\sigma \hat{\mathcal{D}f}\left(\sigma.\phi\right)
\end{split}
\end{equation}
The result follows from the Fourier slice theorem.
\end{proof}
\end{lemma}

In \cite[page 92]{nat1} Natterer explains why it is reasonable to consider picture densities as functions $f$ of compact support in $H^{\alpha}\left(\mathbb{R}^n\right)$ with $\alpha<1/2$. He then gives a bound for the least squares error in his reconstruction from plane integral data in terms of $\rho$, where $\|f\|_{H^{\alpha}}\leq \rho$. With this in mind we will show that the map $f\to \tilde{f}$ is bounded and has a bounded inverse from $H^{\alpha}\to H^{\alpha}$ for $0<\alpha<1$. First, from \cite[page 204]{nat1}, we have the lemma:

\begin{lemma}
\label{lemma2}
Let $\chi\in C_0^{\infty}\left(\mathbb{R}^n\right)$ and let $f\in H^{\alpha}\left(\mathbb{R}^n\right)$. Then the map $f\to \chi f$ is bounded in $H^{\alpha}$ for any $\alpha \in \mathbb{R}$.
\end{lemma}

Now we have our result:

\begin{lemma}
Let $D_r$ be as defined in section \ref{sec1}. Let $f\in C_0^{\infty}\left(D_r\right)$ for some $r>1$ and let $\tilde{f}$ be defined as in Theorem \ref{th1}. Then there exist constants $c\left(\alpha\right)$ and $C\left(\alpha\right)$ such that:
\begin{equation}c\left(\alpha\right)\|\tilde{f}\|_{{H}^{\alpha}}\leq \|f\|_{{H}^{\alpha}}\leq C\left(\alpha\right)\|\tilde{f}\|_{{H}^{\alpha}}\end{equation}
for any $0<\alpha<1$.

\begin{proof}
Let $\chi_{D_r}\in C_0^{\infty}\left(\mathbb{R}^2\right)$ be $1$ on $D_r$ and let $\chi=\chi_{D_r}|x|^4$. Then by Lemma \ref{lemma2}, we have:
\begin{equation}
\begin{split}
c_1\left(\alpha\right)\|f\|^2_{{H}^{\alpha}}&\geq \|\chi_{D_r}|x|^4 f\|^2_{{H}^{\alpha}}\\
&= \||x|^4 f\|^2_{{H}^{\alpha}}\\
&={\|f\|^2_{L^2}}+\iint_{\mathbb{R}^2 \times \mathbb{R}^2}\frac{||x|^4f\left(x\right)-|y|^4f\left(y\right)|^2}{|x-y|^{2+2\alpha}}\mathrm{d}x\mathrm{d}y\\
&\geq {\left(1/r^4\right)\|\tilde{f}\|^2_{L^2}}+\left(1/r^2\right)^{2\alpha-2}\iint_{\mathbb{R}^2 \times \mathbb{R}^2}\frac{|\tilde{f}\left(x\right)-\tilde{f}\left(y\right)|^2}{|x-y|^{2+2\alpha}}\mathrm{d}x\mathrm{d}y\\
&\geq c_2\left(\alpha\right)\|\tilde{f}\|^2_{{H}^{\alpha}}
\end{split}
\end{equation}
for $0<\alpha<1$, This proves the left hand inequality. The right hand inequality can be proven in a similar way.
\end{proof}
\end{lemma}

Before proving the main theorem of this section we state the interpolation inequality for Sobolev spaces on $\mathbb{R}^n$ \cite[page 203]{nat1}.

\begin{lemma}
Let $f\in H^{\gamma}\left(\mathbb{R}^n\right)$. Then we have:
\begin{equation}\|f\|_{H^{\gamma}\left(\mathbb{R}^n\right)}\leq \|f\|^{\frac{\beta-\gamma}{\beta-\alpha}}_{H^{\alpha}\left(\mathbb{R}^n\right)}\|f\|^{\frac{\gamma-\alpha}{\beta-\alpha}}_{H^{\beta}\left(\mathbb{R}^n\right)}\end{equation}
for any $\alpha < \gamma<\beta$.
\end{lemma}

Now we have our main theorem for this section:

\begin{theorem}
\label{the2}
Let $f\in C_0^{\infty}\left(D_r\right)$ for some $r>1$. Then we have:
\begin{equation}\|f\|_{L^2\left(\mathbb{R}^2\right)}\leq c\left(\beta\right)\rho^{\frac{3/2}{\beta+3/2}} \|\mathcal{D}f\|^{\frac{\beta}{\beta+3/2}}_{L^2\left([-1,1]\times S^1\right)}\end{equation}
for any $0<\beta<1$ with $\|f\|_{H^{\beta}}\leq \rho$.

\begin{proof}
Let $\tilde{f}$ be defined as in Theorem \ref{th1} and let $\Phi=\left(\cos \phi,\sin\phi\right)$. Then we have:
\begin{equation}
\begin{split}
2\|\mathcal{D}f\|^{2}_{L^2\left([-1,1]\times S^1\right)}&=2\int_{S^1}\int_{-1}^{1}|\mathcal{D}f\left(p,\phi \right)|^2\mathrm{d}p\mathrm{d}\phi\\
&\geq \int_{S^1}\int_{-1}^{1}|\mathcal{D}f\left(p,\phi \right)+\mathcal{D}f\left(-p,\phi \right)|^2\mathrm{d}p\mathrm{d}\phi\\
&= \int_{S^1}\int_{-\infty}^{\infty}|\hat{\mathcal{D}f}\left(\sigma,\phi \right)+\hat{\mathcal{D}f}\left(-\sigma,\phi \right)|^2\mathrm{d}\sigma\mathrm{d}\phi\\
&= 2\pi\int_{S^1}\int_{-\infty}^{\infty}\frac{1}{\sigma^2}|\hat{\tilde{f}}\left(\sigma\Phi\right)-\hat{\tilde{f}}\left(-\sigma\Phi\right)|^2\mathrm{d}\sigma\mathrm{d}\phi\\
&\geq 2\pi\int_{S^1}\int_{-\infty}^{\infty}\left(1+\sigma^2\right)^{-1}|\hat{\tilde{f}}\left(\sigma\Phi\right)-\hat{\tilde{f}}\left(-\sigma\Phi\right)|^2\mathrm{d}\sigma\mathrm{d}\phi\\
&= 4\pi\int_{S^1}\int_{0}^{\infty}\left(1+\sigma^2\right)^{-1}|\hat{\tilde{f}}\left(\sigma\Phi\right)-\hat{\tilde{f}}\left(-\sigma\Phi\right)|^2\mathrm{d}\sigma\mathrm{d}\phi
\end{split}
\end{equation}
After making the substitution $\xi=\sigma \Phi$, we have:
\begin{equation}
\begin{split}
\frac{1}{2\pi}\|\mathcal{D}f\|^{2}_{L^2\left([-1,1]\times S^1\right)}&\geq \int_{\mathbb{R}^2}|\xi|^{-1}\left(1+|\xi|^2\right)^{-1}|\hat{\tilde{f}}\left(\xi\right)-\hat{\tilde{f}}\left(-\xi\right)|^2\mathrm{d}\xi\\
&\geq \|\tilde{f}-\tilde{f}^{-}\|^2_{H^{-\frac{3}{2}}\left(\mathbb{R}^2\right)}
\end{split}
\end{equation}
where $\tilde{f}^{-}\left(x\right)=\tilde{f}\left(-x\right)$. Applying the interpolation inequality with $\alpha=-3/2$ and $\gamma=0$, yields:
\begin{equation}
\begin{split}
2\|f\|_{L^2\left(\mathbb{R}^2\right)}&\leq 2\|\tilde{f}\|_{L^2\left(\mathbb{R}^2\right)}\\
&=\|\tilde{f}-\tilde{f}^{-}\|_{L^2\left(\mathbb{R}^2\right)}\ \ \ \text{since}\ \ \ f\in C_0^{\infty}\left(D_r\right)\\
&\leq \|\tilde{f}-\tilde{f}^{-}\|^{\frac{\beta}{\beta+3/2}}_{H^{-\frac{3}{2}}\left(\mathbb{R}^2\right)}\|\tilde{f}-\tilde{f}^{-}\|^{\frac{3/2}{\beta+3/2}}_{H^{\beta}\left(\mathbb{R}^2\right)}\\
&\leq \left(2\pi\right)^{-1/2}2^{\frac{3/2}{\beta+3/2}}\|\mathcal{D}f\|^{\frac{\beta}{\beta+3/2}}_{L^2\left([-1,1]\times S^1\right)}\|\tilde{f}\|^{\frac{3/2}{\beta+3/2}}_{H^{\beta}\left(\mathbb{R}^2\right)}\\
&\leq c\left(\beta\right)\|\mathcal{D}f\|^{\frac{\beta}{\beta+3/2}}_{L^2\left([-1,1]\times S^1\right)}\|f\|^{\frac{3/2}{\beta+3/2}}_{H^{\beta}\left(\mathbb{R}^2\right)}\\
&\leq c\left(\beta\right)\rho^{\frac{3/2}{\beta+3/2}} \|\mathcal{D}f\|^{\frac{\beta}{\beta+3/2}}_{L^2\left([-1,1]\times S^1\right)}
\end{split}
\end{equation}
for any $0<\beta<1$ with $\|f\|_{H^{\beta}}\leq \rho$.
\end{proof}
\end{theorem}

\begin{corollary}
Let $f\in C_0^{\infty}\left(D_r\right)$ for some $r>1$ and let $g=\mathcal{D}f$. Let $g^{\epsilon}\in L^2\left([-1,1]\times S^1\right)$ be such that $\|g^{\epsilon}-g\|_{L^2\left([-1,1]\times S^1\right)}<\epsilon$. Then for any $f_1,f_2 \in C_0^{\infty}\left(D_r\right)$ which satisfy $\|\mathcal{D}f-g^{\epsilon}\|_{L^2\left([-1,1]\times S^1\right)}<\epsilon$, we have:
\begin{equation}\|f_1-f_2\|_{L^2\left(\mathbb{R}^2\right)} \leq c\left(\beta\right) \rho^{\frac{3/2}{\beta+3/2}}\epsilon^{\frac{\beta}{\beta+3/2}}\end{equation}
for any $0<\beta<1$ with $\|f_i\|_{H^{\beta}}\leq \rho$ for $i=1,\ 2$.
\end{corollary}

We can interpret this last corollary to mean that given some erroneous data $g^{\epsilon}$ which differs in the least squares sense from $\mathcal{D}f$ absolutely by $\epsilon$, the least squares error in our solution is bounded above by $c\left(\beta\right) \rho^{\frac{3/2}{\beta+3/2}}\epsilon^{\frac{\beta}{\beta+3/2}}$ for some constant $c\left(\beta\right)$ with the a-priori knowledge that $\|f\|_{H^{\beta}}\leq \rho$. 

In \cite{nat1} Natterer uses the value $\beta/\left(\alpha+\beta\right)$ as a measure for the ill posedness of his problem and gives his criteria for a linear inverse problem to be modestly, mildly or severely ill posed. If we set $\beta$ close to $1/2$, then based on these criteria the above arguments would suggest that our problem is mildly ill posed, but more ill posed than the inverse Radon transform, which we would expect given that the disc transform $\mathcal{D}$ is a degree smoother than $R$. 
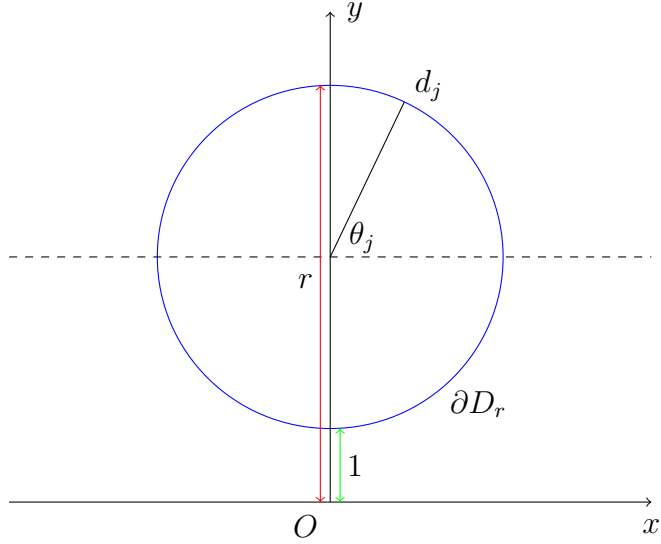
\begin{figure}[!h]
\centering
\begin{tikzpicture}[scale=3.25]
\draw [blue] (0,0) circle [radius=0.7];
\draw [->, thin] (-1.3,-1)--(1.3,-1);
\node at (-0.1,-1.1) {$O$};
\draw [->,thin] (0,-1)--(0,1);
\node at (0.1,1) {$y$};
\node at (1.3,-1.1) {$x$};
\draw [dashed, thin] (-1.3,0)--(1.3,0);
\draw [<->,thin, green] (0.04,-1)--(0.04,-0.7);
\draw [<->,thin, red] (-0.04,-1)--(-0.04,0.7);
\node at (-0.1,-0.1) {$r$};
\node at (0.1,-0.85) {$1$};
\node at (0.6,-0.6) {$\partial D_r$};
\draw [thin] (0,0)--(0.3,0.632);
\node at (0.13,0.08) {$\theta_j$};
\node at (0.4,0.7) {$d_j$};
\end{tikzpicture}
\caption{A representation of the detector ring $\partial D_r$ in the proposed coordinate system. The polar angle $\theta_j \in [0,2\pi]$ determines the  detector position $d_j$.}
\label{fig0}
\end{figure}

Another source of error in our solution can be due to limited sampling of the data. In practice the number of detectors will be finite. Let us parameterize the set of points on the detector ring $\partial D_r=\{\left(x,y\right)\in \mathbb{R}^2 : x^2+\left(y-\left(r+1\right)/2\right)^2=\left(r-1\right)^2/4\}$ in terms of a polar angle $\theta$, and let the finite set of polar angles $\Theta=\{\theta_1,\ldots,\theta_n\}$ determine a finite set of detector positions $\{d_1,\ldots,d_n\}\in \partial D_r$. See figure \ref{fig0}. Then for every $\phi\in [0,2\pi]$ we can sample $\mathcal{D}f(p,\phi)$ for:
\begin{equation}
\label{sample}
p=p_j=\frac{r\cos\theta_j \sin \phi +\left(1+r\sin\theta_j\right)\sin \phi}{r^2+1+2r\sin\theta_j},\ \ \  \ \ 1\leq j\leq n
\end{equation}
where $p_j$ is such that $\{\frac{1}{2}\left((r-1)\cos(\theta_j), (1+\sin\theta_j)r+1-\sin \theta_j\right)\}\subset \partial D_{\frac{1}{p_j},\phi}\cap \partial D_r$ for $1\leq j\leq n$.

From \cite[pages 204 and 42]{nat1} we have:

\begin{lemma}
\label{lemma10}
Let $\Omega\subset \mathbb{R}^n$ be bounded and sufficiently regular. For $h>0$ let $\Omega_k$ be a finite subset of $\Omega$ such that $d\left(\Omega,\Omega_k\right)\leq h$, where $d$ is the is the Hausdorff distance metric between sets. Let $\alpha>n/2$ where $\alpha=m+\sigma$ for some integer $m$ and $0<\sigma <1$, and for $f\in H^{\alpha}\left(\Omega\right)$  define the seminorm:
\begin{equation}|f|^2_{H^{\alpha}\left(\Omega\right)}=\sum_{|k|=m}\iint_{\Omega\times \Omega}\frac{|D^kf\left(x\right)-D^kf\left(y\right)|^2}{|x-y|^{n+2\sigma}}\mathrm{d}x\mathrm{d}y\end{equation}
Then there is a constant $c$ such that:
\begin{equation}\|f\|_{L^2\left(\Omega\right)}\leq ch^{\alpha}|f|_{H^{\alpha}\left(\Omega\right)}\end{equation}
for every $f\in H^{\alpha}\left(\Omega\right)$ which is zero on $\Omega_k$.
\end{lemma}

\begin{theorem}
Let $\Omega_n$ be the unit ball in $\mathbb{R}^n$. For every $\alpha$ there exist positive constants $c\left(\alpha,n\right)$ and $C\left(\alpha,n\right)$ such that for $f\in C_0^{\infty}\left(\Omega_n\right)$
\begin{equation}c\left(\alpha,n\right)\|f\|_{H^{\alpha}\left(\Omega_n\right)}\leq \|Rf\|_{H^{\alpha+\left(n-1\right)/2}\left(Z\right)}\leq C\left(\alpha,n\right)\|f\|_{H^{\alpha}\left(\Omega_n\right)}\end{equation}
\end{theorem}

From these we have the theorem:

\begin{theorem}
For each $\phi\in [0,2\pi]$ let $I_{\phi}\subset [-1,1]$ be a finite subset of the unit interval. Let:
\begin{equation}h=\sup_{\phi} d\left(I_{\phi},[-1,1]\right)\end{equation}
where $d$ is the Hausdorff distance metric. Let $f\in C_0^{\infty}\left(D_r\right)$ and let $\|f\|_{H^{\alpha}}<\rho$ with $0<\alpha<1$. If $\mathcal{D}f_{\phi}$ is zero on $I_{\phi}$ for every $\phi \in[0,2\pi]$, then there exists a constant $c\left(\alpha\right)$ such that:
\begin{equation}\|f\|_{L^2\left(\mathbb{R}^2\right)}\leq c\left(\alpha\right)h^{\alpha}\rho\end{equation}

\begin{proof}
Let $\tilde{f}$ be defined as in Theorem \ref{th1} and let $|\cdot|_{H^{\alpha}}$ be the seminorm defined in Lemma \ref{lemma10}. Let $\alpha+3/2=m+\sigma$ for some integer $m$ and $0<\sigma<1$. Then we have:
\begin{equation}
\begin{split}
\|\mathcal{D}f \|^2_{L^2\left([-1,1]\times S^1\right)}&=\int_{S^1}\|\mathcal{D}f_{\phi}\|^2_{L^2\left([-1,1]\right)} \mathrm{d}\phi\\
&\leq c^2 h^{2\alpha+3}\int_{S^1}|\mathcal{D}f_{\phi}|^2_{H^{\alpha+3/2}\left([-1,1]\right)} \mathrm{d}\phi\\
&\leq c^2 h^{2\alpha+3}\int_{S^1}\iint_{[-1,1]\times [-1,1]}\frac{|\frac{\partial^m}{\partial p^m} \mathcal{D}f_{\phi}-\frac{\partial^m}{\partial p^m} \mathcal{D}f_{\phi}|^2}{|x-y|^{n+2\sigma}}\mathrm{d}x\mathrm{d}y \mathrm{d}\phi\\
&= c^2 h^{2\alpha+3}\int_{S^1}\iint_{[-1,1]\times [-1,1]}\frac{|\frac{\partial^{m-1}}{\partial p^{m-1}} R\tilde{f}_{\phi}-\frac{\partial^{m-1}}{\partial p^{m-1}} R\tilde{f}_{\phi}|^2}{|x-y|^{n+2\sigma}}\mathrm{d}x\mathrm{d}y \mathrm{d}\phi\\
&= c^2 h^{2\alpha+3}\int_{S^1}|R\tilde{f}_{\phi}|^2_{H^{\alpha+1/2}\left([-1,1]\right)} \mathrm{d}\phi\\
&\leq c^2 h^{2\alpha+3}\int_{S^1}\|R\tilde{f}_{\phi}\|^2_{H^{\alpha+1/2}\left([-1,1]\right)} \mathrm{d}\phi\\
&= c^2 h^{2\alpha+3}\|R\tilde{f}\|^2_{H^{\alpha+1/2}\left(Z\right)}\\
&\leq c_1\left(\alpha\right)^2 h^{2\alpha+3}\|\tilde{f}\|^2_{H^{\alpha}\left(\mathbb{R}^2\right)}\\
&\leq c_2\left(\alpha\right)^2 h^{2\alpha+3}\|f\|^2_{H^{\alpha}\left(\mathbb{R}^2\right)}\\
&\leq c_2\left(\alpha\right)^2 h^{2\alpha+3}\rho^2\\
\end{split}
\end{equation}
for $0<\alpha<1$ with $\|f\|_{H^{\alpha}}\leq \rho$. Applying Theorem \ref{the2}, we have:
\begin{equation}
\begin{split}
\|f\|_{L^2\left(\mathbb{R}^2\right)}&\leq c_3\left(\alpha\right)\rho^{\frac{3/2}{\alpha+3/2}} \|\mathcal{D}f\|^{\frac{\alpha}{\alpha+3/2}}_{L^2\left([-1,1]\times S^1\right)}\\
&\leq c\left(\alpha\right)h^{\alpha}\rho
\end{split}
\end{equation}
which completes the proof.
\end{proof}
\end{theorem}

This last result tells us that given a finite set of detectors with a disc diameter sampling determined by equation (\ref{sample}) and with $h$ being a measure of the uniformity of the sample, the least squares error in our solution is bounded above by $c\left(\alpha\right)h^{\alpha}\rho$ with the a-priori knowledge that $\|f\|_{H^{\alpha}}\leq \rho$ for some $0<\alpha<1$.





\section{The physical model}
\label{phys}
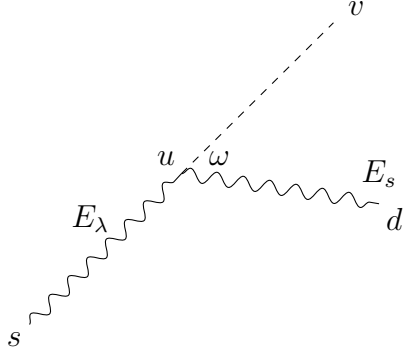
\begin{figure}[!h]
\centering
\begin{tikzpicture}[scale=2]
\draw [snake it, thin] (-1,-1)--(0,0);
\draw [dashed, thin] (0,0)--(1,1);
\draw [snake it, thin] (0,0)--(1.3,-0.2);
\node at (0.25,0.1) {$\omega$};
\node at (-1.1,-1.1) {$s$};
\node at (1.4,-0.3) {$d$};
\node at (-0.6,-0.3) {$E_{\lambda}$};
\node at (1.3,0) {$E_s$};
\node at (-0.1,0.1) {$u$};
\node at (1.15,1.1) {$v$};
\end{tikzpicture}
\caption{A scattering event with initial photon energy $E_{\lambda}$ from a source $s$ scattered to $d$ with energy $E_s$. The dashed line displays the original path of the photon to a detector $v$.} 
\label{figurep}
\end{figure}

\begin{figure}[!h]
\centering
\includegraphics[scale=0.5]{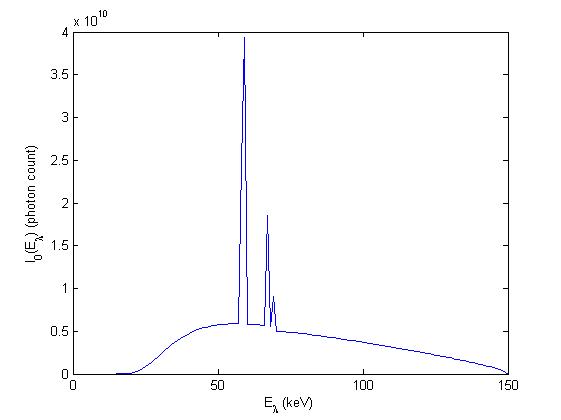}
\caption{A typical Tungsten target spectrum with a 150kV accelerating voltage.}
\label{figspec}
\end{figure}

In this section we present an accurate physical model and  a possible approximate model which allows for the proposed reconstruction method. We consider an intensity of photons scattering from a point $u$ as illustrated in figure \ref{figurep}.
The intensity of photons scattered from $u$ to $d$ with energy $E_s$ is:
\begin{equation}
\label{equ17}
\begin{split}
I\left(u,d,E_s\right)=I_0\left(E_{\lambda}\right)
\exp\left(-\int_{l_{1}} \mu_{E_{\lambda}}\right)&n_e\left(u\right)\mathrm{d}V\\
&\times \frac{\mathrm{d}\sigma}{\mathrm{d}\Omega}\left(E_s,\omega\right)S\left(q\right) 
\exp\left(-\int_{l_2} \mu_{E_s}\right)\mathrm{d}\Omega_{u,d}
\end{split}
\end{equation}
where $I_0$ is the initial intensity, which depends on the energy $E_{\lambda}$ (see figure \ref{figspec} for an example polychromatic spectrum). $\mu_{E}$ is the linear attenuation coefficient, which is dependant on the energy $E$ and the atomic number of the target material. Here $n_e\left(u\right)\mathrm{d}V$ is the number of electrons in a volume $\mathrm{d}V$ around the scattering point $u$. So $n_e$ (number of electrons per unit volume) is the quantity to be reconstructed. $l_1$ and $l_2$ are the line segments connecting $s$ to $u$ and $u$ to $d$ respectively.

The Klein-Nishina differential cross section $\mathrm{d}\sigma/\mathrm{d}\Omega$, is defined by:\begin{equation}\frac{\mathrm{d}\sigma}{\mathrm{d}\Omega}\left(E_s,\omega\right)=\frac{r_0^{2}}{2}{\left(\frac{E_s}{E_{\lambda}}\right)}^{2}\left(\frac{E_s}{E_{\lambda}}+\frac{E_{\lambda}}{E_s}-1+\cos^{2} \omega\right)\end{equation}
where $r_0$ is the classical electron radius. This predicts the scattering distribution for a photon off a free electron at rest. Given that the atomic electrons typically are neither free nor at rest, a correction factor is included, namely the incoherent scattering function $S\left(q\right)$. Here $q=\frac{E_{\lambda}}{hc}\sin \left(\omega /2\right)$ is the momentum transferred by a photon with initial energy:
\begin{equation}E_{\lambda}=\frac{E_s}{1-\left(E_s/E_0\right)\left(1-\cos\omega\right)}\end{equation}
scattering at an angle $\omega$, where $h$ is Planck's constant and $c$ is the speed of light. The scattering function $S$ also depends on the atomic number $Z$, so we set $Z=Z_{\text{avg}}$ to some average atomic number as an approximation. 


For $Z_{\text{avg}}=45$ (Rhodium) we have the expression:
\begin{equation}
\label{equ16}
S\left(q\right)=1-\frac{1.023}{\left(1+0.458q\right)^{2.509}}
\end{equation}
To acquire equation (\ref{equ16}) we have extended the least squares fit given in \cite{fit} to the values of $S\left(q\right)$ given in \cite{hub}.

The solid angle subtended by $u$ and $d$ is defined:
\begin{equation}\mathrm{d}\Omega_{u,d}=\frac{A}{4\pi}\frac{\bold{r}\cdot  \bold{n}}{|\bold{r}|^3}\end{equation}
where $\bold{r}=d-u$, $A$ is the detector area and $\bold{n}$ is the unit vector normal to the detector surface. 

Given our machine geometry and proposed reconstruction method, it is difficult to include the more accurate model stated above as an additional weighting to our integral equations (as in done in \cite{norton} for example) while allowing for the same inversion formulae. So we average equation (\ref{equ17}) over the scattering region $\text{R}_{p,\phi}=D_{\frac{1}{p},\phi}\cap D_r$, for each $p$ and $\phi$. Here $D_{p,\phi}$ and $D_r$ are as defined in section \ref{sec1}, where $r$ is fixed depending on the machine specifications.

Let $I\left(u,d,E_s\right)=I\left(u;p,\phi\right)=P\left(u;p,\phi\right) n_e\left(u\right) \mathrm{d}V$. Here $P=P\left(u;p,\phi\right)$ depends on the scattering point $u$ and $p$ and $\phi$ as defined in section \ref{sec1}. When $\text{R}_{p,\phi}\neq \emptyset$, $p$ and $\phi$ determine the detector position $d$ and the measured energy $E_s$. 

We have:
\begin{equation}
\label{app2}
P_{\text{avg}}\left(p,\phi\right)=\frac{1}{{A}\left(\text{R}_{p,\phi}\right)}\iint_{\text{R}_{p.\phi}} P\left(u;p,\phi\right) \mathrm{d}u
\end{equation}
which gives the average of $P$ over $\text{R}_{p,\phi}$. Here ${A}\left(\text{R}_{p,\phi}\right)$ denotes the area of $\text{R}_{p,\phi}$. Let $f:\mathbb{R}^2\to \mathbb{R}$ be an example density with support contained in $D_r$, and let:
\begin{equation}
I_1\left(C,p,\phi\right)=sC\iint_{\text{R}_{p.\phi}} P\left(u;p,\phi\right) \mathrm{d}u
\end{equation}
be the scattered intensity measured for a constant density $C$ over $\text{R}_{p,\phi}$, where $s$ is the (constant) slice thickness. Then if we assume that the scattering probability is constant and equal to $P_{\text{avg}}(p,\phi)$ throughout each scattering region $\text{R}_{p.\phi}$, the absolute error in our approximation would satisfy:
\begin{equation}
\left|I_m(p,\phi)-sP_{\text{avg}}\left(p,\phi\right)\mathcal{D}_1f(p,\phi)\right|\leq \left|I_1\left(\max_{x\in \text{R}_{p,\phi}}f(x),p,\phi\right)-I_1\left(\min_{x\in \text{R}_{p,\phi}}f(x),p,\phi\right)\right|
\end{equation}
for all $(p,\phi)\in Z^{+}$. Here $I_m : Z^{+}\to \mathbb{R}$ is the intensity of photons we measure. So provided that the range of the density values is small over the majority of scattering regions considered, the averaged model given above will have a similar level of accuracy to the more precise model given in equation (\ref{equ17}).

If the linear attenuation coefficient $\mu$ is known a-priori, then the exponential terms of equation (\ref{equ17}) may be included in $P$. Otherwise we may approximate:
\begin{equation}\exp\left(-\int_{l_{1}} \mu_{E_{\lambda}}\right)\exp\left(-\int_{l_2} \mu_E\right)\approx\exp\left(-\int_{l_v}\mu_{E_{\lambda}}\right)\end{equation}
where $l_v$ is the line segment from $s$ to the detector in the forward direction $v$ (see figure \ref{figurep}). This is the approximation made  in \cite{Nic}.  By the Beer-Lambert law, we have:
\begin{equation}\frac{I_v\left(E_{\lambda}\right)}{I_0\left(E_{\lambda}\right)}=\exp\left(-\int_{l_{v}}\mu_{E_{\lambda}}\right)\end{equation}
where $I_v\left(E_{\lambda}\right)$ is the recorded straight through intensity. 

To account for the physical modelling, we would divide the data by $sP_{\text{avg}}$ to calculate approximate values for $\mathcal{D}_1f$ and hence for $\mathcal{D}f$. 


\subsection{Determining the atomic number}
\label{sec4}
With the proposed machine configuration, we can show that the data collected in the light field determines the linear attenuation coefficient $\mu$ uniquely (this is the standard 2D reconstruction problem). With the additional information provided by our theory, we show under the right assumptions that the atomic number of the target is determined uniquely by the full data (light plus dark field).

\begin{figure}[!h]
\centering
\includegraphics[scale=0.6]{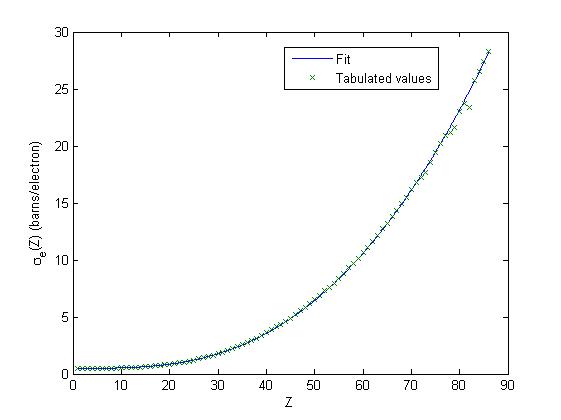}
\caption{We have presented our fit for $\sigma_e$ for $E=100$keV up to $Z=86$ where a sudden dip in the $\sigma_e$ values occurs. The tabulated values of $\sigma_e$ given in \cite{vieg} are shown alongside the fitted curve.}
\label{F6}
\end{figure}

The electron density $n_e$ and the linear attenuation coefficient $\mu$ are related via the formula:
\begin{equation}
\label{equ30}
\mu\left(E,Z\right)=n_e \sigma_e\left(E,Z\right)
\end{equation}
where $\sigma_e$ is the total cross section per electron. The cross section $\sigma_e$ is continuous and monotone increasing as a function of $Z$ on $[1,Z_{\text{max}}]$, where $Z_{\text{max}}\propto\sqrt{E}$.  For example, when $E=100$keV we can calculate $Z_{\text{max}}= 86$. With this in mind, we fit a smooth curve to known values of $\sigma_e$ given in the atomic data tables \cite{vieg}. This allows us to calculate values of $\sigma_e$ for non integer $Z$. See figure \ref{F6}.

The formula for the fit presented for a fixed energy $E=100$keV is:
\begin{equation}
\sigma_e\left(Z\right)=\sigma_p\left(Z\right)+\sigma_s\left(Z\right)=\left(1.51\times10^{-6} Z^{4.72-0.22\log Z}\right) +\left(0.49+7.90\times 10^{-4}\left(1-Z^{-0.50}\right)Z^{1.57}\right)
\end{equation}
This was obtained via a combination of the formula for $\sigma_s$ (the total scattering cross section) presented by Jackson and Hawkes in \cite{DJ} and the suggested fit for $\sigma_p$ (the photoelectric cross section) given in \cite{vieg}. Fits for energies other than $E=100$keV are also possible via the same fitting method.

From our theory we know that the data determines $\mu_E$ and $n_e$ uniquely, where $E\leq E_{\text{max}}$. If we assume that atomic numbers $Z\geq Z_{\text{max}}$ are not present in the target material, then it is clear from the above arguments that the atomic number of the target is uniquely determined. Without this assumption the atomic number would be limited to a range of values. If we reconstruct both $\mu_E$ and $n_e$ for a suitably high energy $E$, we can then calculate values for $Z$ from our curve fit for $\sigma_e$. We will test this additional method in our results also.


\section{Results}
\label{res}
\begin{figure}[!h]
\includegraphics[scale=0.45]{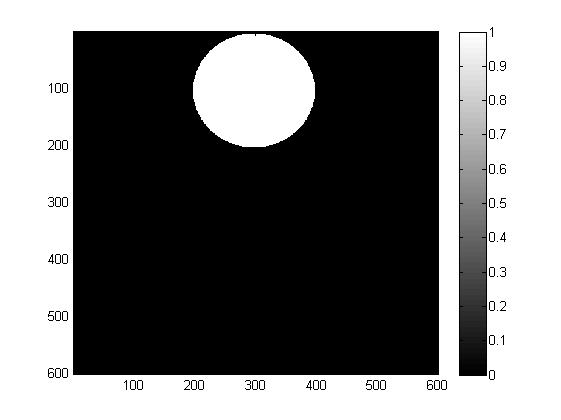}
\includegraphics[scale=0.45]{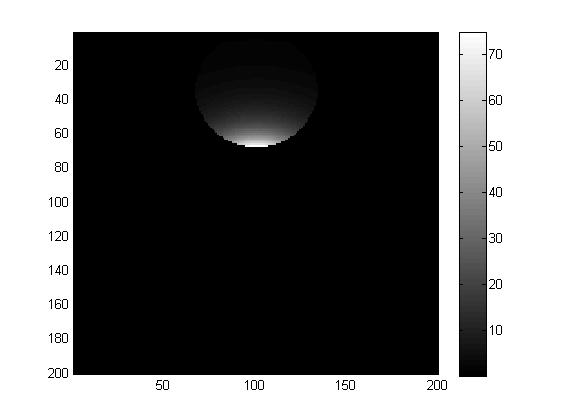}
\caption{A water bottle cross section $f$ is simulated as a circular region of uniform density $1$ on the left. The function $\tilde{f}$ as defined in Theorem \ref{th1} is shown on the right.}
\label{figure5}
\end{figure}

\begin{figure}[!h]
\includegraphics[scale=0.45]{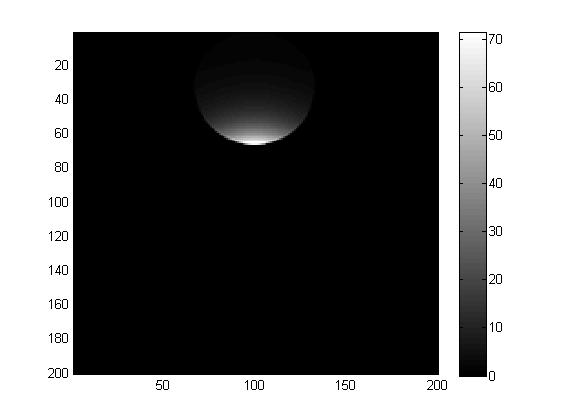}
\includegraphics[scale=0.45]{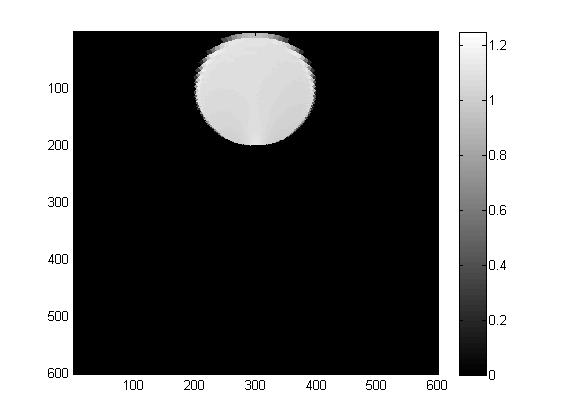}
\caption{A reconstruction of $\tilde{f}$ in the absence of added noise is shown on the left. We have backprojected from 180 views with the default Ram-Lak cropped filter. The corresponding pixel values of $f$ are presented on the right. Both $f$ and $\tilde{f}$ are set to $0$ outside of their support. 
}
\label{figure7}
\end{figure}

To test our reconstruction methods, let us consider the water bottle cross section $f$ and the corresponding function $\tilde{f}$ represented in figure \ref{figure5}. We calculate values of $\mathcal{D}f$ for $p$ in the range $[0,1]$ and for $\phi \in \{\frac{\pi}{180},\cdots,2\pi\}$. These were calculated using the exact formula for the area of intersection of two discs. We approximate the derivative of $\mathcal{D}f$ with respect to $p$ as the finite difference:
\begin{equation}
\label{app1}
\frac{\partial}{\partial p}\mathcal{D}f\left(p,\phi\right)=\frac{\mathcal{D}f\left(p+h,\phi\right)-\mathcal{D}f\left(p,\phi\right)}{h}
\end{equation}
for a chosen step size $h$. 
To reconstruct $\tilde{f}$ we apply the Matlab function ``iradon", which filters (choosing from a selection of filters pre-coded by Matlab) and backprojects the projection data $R\tilde{f}=-\frac{\partial}{\partial p}\mathcal{D}f$ to recover $\tilde{f}$. We then make the necessary change in coordinates to produce our density image $f$. In the absence of noise we find our results to be satisfactory. See figure \ref{figure7}. 

\begin{figure}[!h]
\includegraphics[scale=0.45]{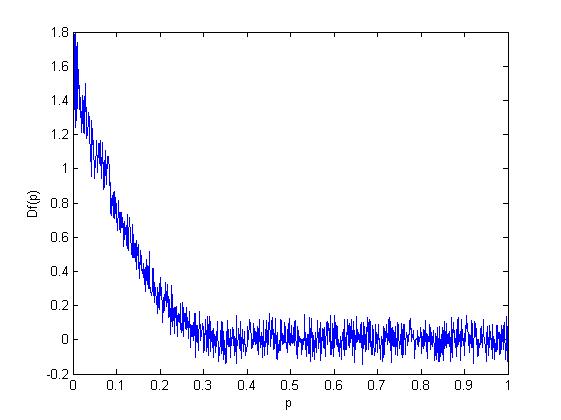}
\includegraphics[scale=0.45]{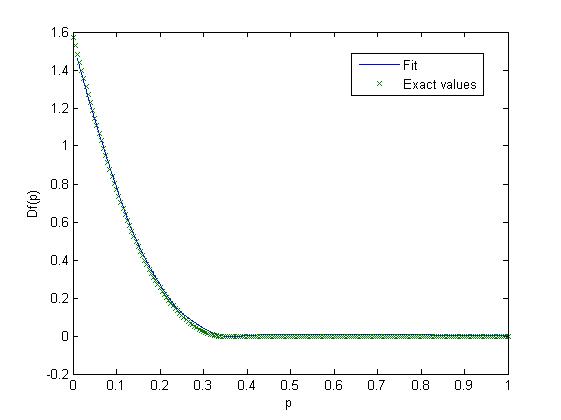}
\caption{On the left we have plotted values of $\mathcal{D}f\left(p,0\right)$ for $p\geq0$ with $10\%$ random noise added. On the right we have applied a simple moving average filter to the simulated data and taken a subsample of the smoothed data before interpolating as specified earlier. The exact values are presented alongside the fitted values in the right hand figure.}
\label{figure9}
\end{figure}

\begin{figure}[!h]
\includegraphics[scale=0.45]{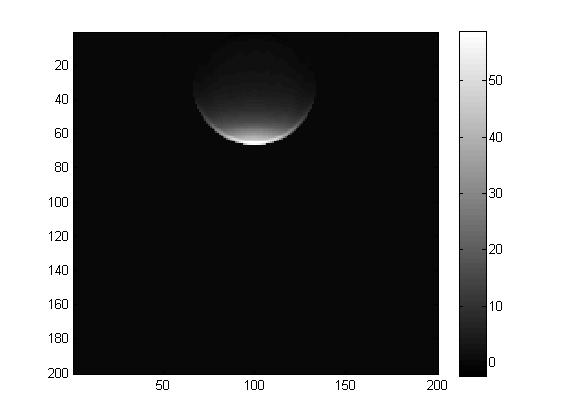}
\includegraphics[scale=0.45]{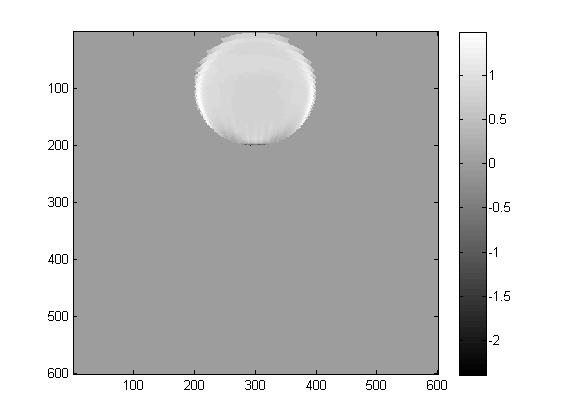}
\caption{On the left we have a reconstruction of $\tilde{f}$ after smoothing with $2\%$ added noise. We have again backprojected from 180 views, although here we have multiplied the standard ramp filter by a Hamming window to reduce the high frequency noise. The corresponding pixel values for $f$ are presented on the right.}
\label{figure11}
\end{figure}

\begin{figure}[!h]
\includegraphics[scale=0.45]{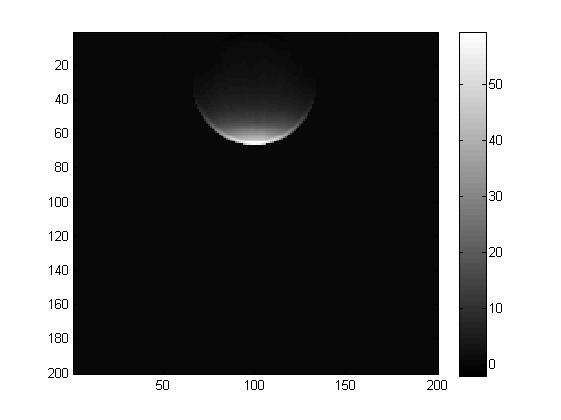}
\includegraphics[scale=0.45]{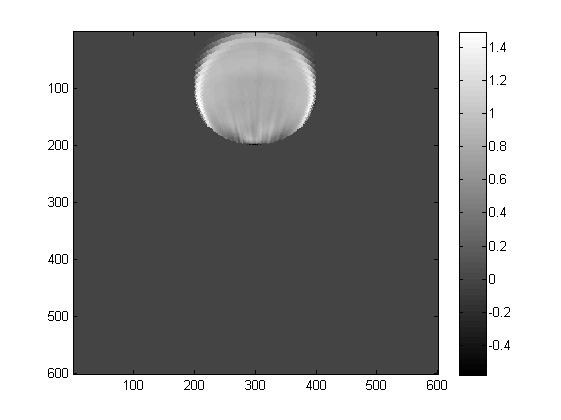}
\caption{On the left, a reconstruction of $\tilde{f}$ after smoothing with $10\%$ added noise. We have multiplied the ramp filter by a Hamming window and backprojected from 180 views. The corresponding pixel values for $f$ are displayed on the right.}
\label{figure13}
\end{figure}

\begin{figure}[!h]
\includegraphics[scale=0.45]{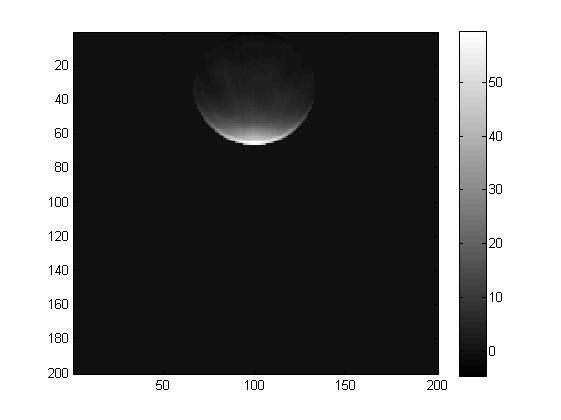}
\includegraphics[scale=0.45]{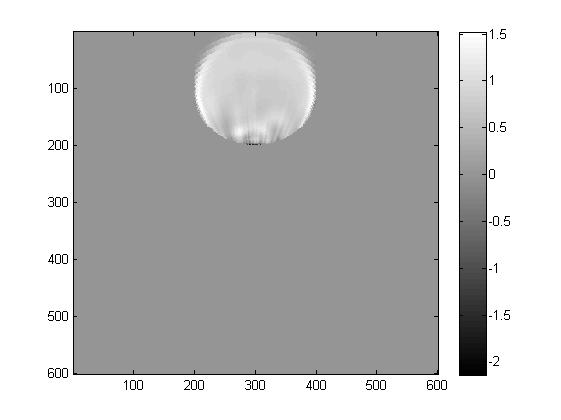}
\caption{A reconstruction of $\tilde{f}$ after smoothing with $50\%$ added noise is shown on the left. We have multiplied the ramp filter by a Hamming window and backprojected from 180 views. The corresponding pixel values for $f$ are displayed on the right.}
\label{figure14}
\end{figure}

Let us now perturb the calculated values of $\mathcal{D}f$ slightly such as to simulate random noise. We multiply each exact value of $\mathcal{D}f$ by a pseudo random number in the range $[1-\frac{\%\text{err}}{100},1+\frac{\%\text{err}}{100}]$ (we use the C++ function ``rand" to generate random numbers), where $\%\text{err}$ is the desired amount of percentage error to be added. In this case, even with a relatively small amount of added noise, the data must be smoothed sufficiently before applying approximation (\ref{app1}). To smooth the data, we apply a simple moving average filter and calculate any intermediate values via a shape preserving cubic interpolation method (``pchip" interpolation in Matlab). We expect this interpolation method to preserve the monotonicity of the data (monotone decreasing) as a function of $p$. To illustrate this technique we refer to figure \ref{figure9}. We have presented our reconstructions after smoothing with $2\%$, $10\%$ and $50\%$ added noise in figures \ref{figure11}, \ref{figure13} and \ref{figure14}. 

\begin{figure}[!h]
\includegraphics[scale=0.45]{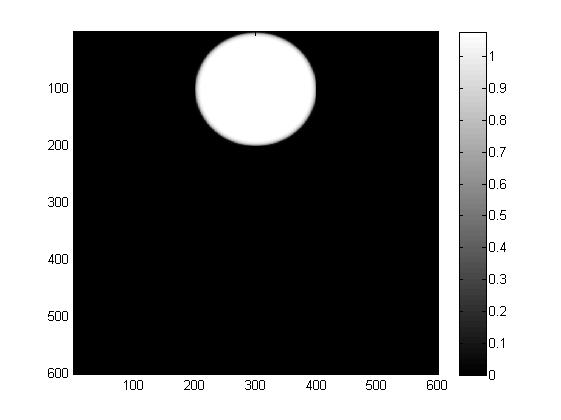}
\includegraphics[scale=0.45]{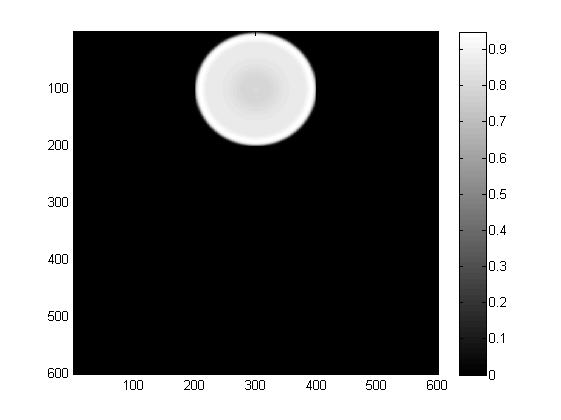}
\caption{On the left, an average reconstruction of $f$ is shown with no noise added to each dataset before reconstruction.  For the right hand image $2\%$ random noise was added to each dataset before reconstruction. In this case $f_{\text{avg}}=0.853$ which gives an atomic number value of $Z=13.3$.}
\label{figure16}
\end{figure}

\begin{figure}[!h]
\includegraphics[scale=0.45]{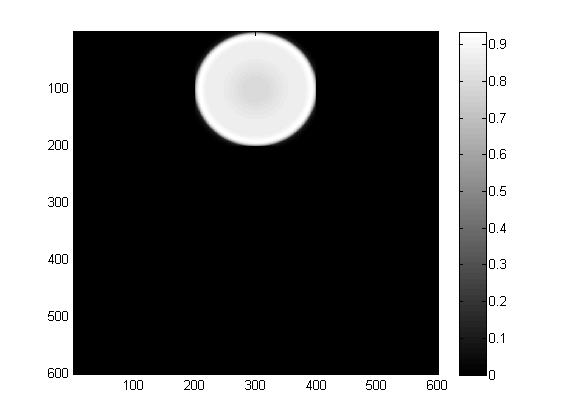}
\includegraphics[scale=0.45]{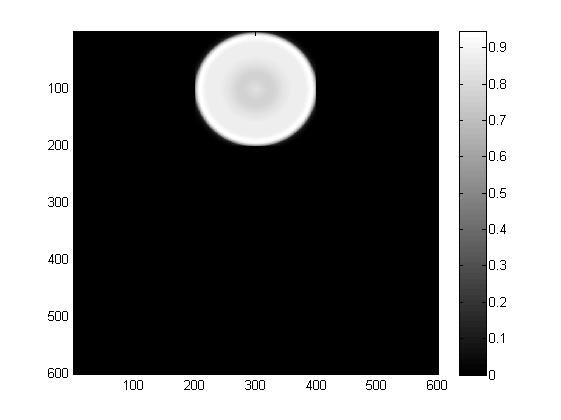}
\caption{On the left, an average reconstruction of $f$ is shown with $10\%$ noise added to each dataset before reconstruction. Here $f_{\text{avg}}=0.865$ which gives an atomic number value of $Z=12.9$. For the right hand image $50\%$ random noise was added to each dataset before reconstruction. In this case $f_{\text{avg}}=0.863$ which gives a reconstructed atomic number value of $Z=13.0$.}
\label{figure15}
\end{figure}

Here, we have reconstructed $f$ from a single view point, using data collected from a single source projection. With the proposed machine configuration however, there are a number of views from which $f$ may be reconstructed. So we take an average over 360 views (for source positions at equal $\pi/180$ intervals over the range $[0,2\pi]$). Our results are presented in figures \ref{figure16} and \ref{figure15}. 
Here we see an improvement in the signal-to-noise-ratio.  The rotational symmetry of $f$ about the centre of the circular region of $f$'s support is also recovered.

Let $f_{\text{avg}}$ be the average of the non zero pixel values shown in the left hand image of figure \ref{figure16} and let $Z=7.420$ be the effective atomic number for water. Then we can calculate $f_{\text{avg}}\approx 1.033$ and using equation (\ref{equ30}) we can calculate the total cross section to be:
\begin{equation}
\sigma_e\left(E,7.420\right)=\frac{\mu(E,7.420)}{1.033}=0.493
\end{equation}
for $E=100$keV assuming no additional error. Based on our curve fit for $\sigma_e(100,Z)$, this would yield a reconstructed atomic number of $Z=0.886$, which differs from the accepted value by $88\%$. For the remaining averaged density reconstructions the $f_{\text{avg}}$ and $Z$ values are given in the figure caption.

\begin{figure}[!h]
\includegraphics[scale=0.45]{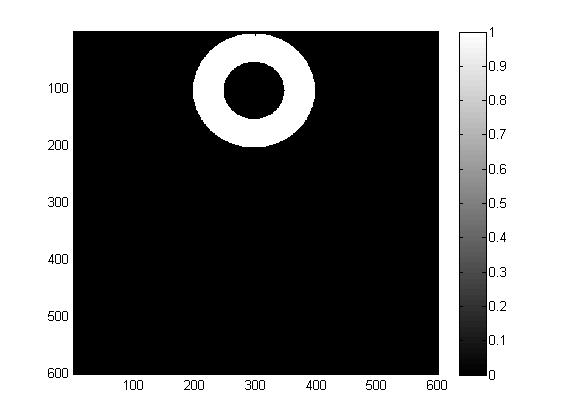}
\includegraphics[scale=0.45]{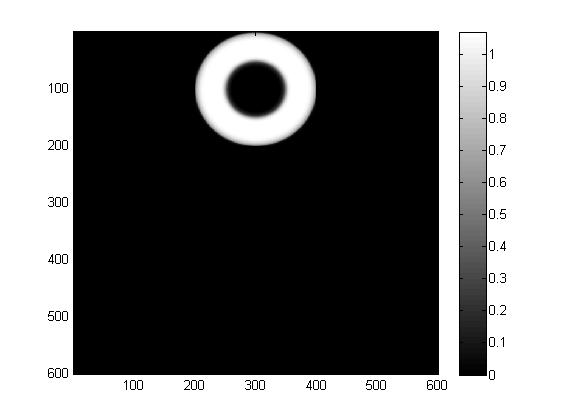}
\caption{On the left we have a simulated hollow tube $f_t$. On the right is an averaged reconstruction of $f_t$ with no noise added to each dataset.}
\label{figure17}
\end{figure}

\begin{figure}[!h]
\includegraphics[scale=0.45]{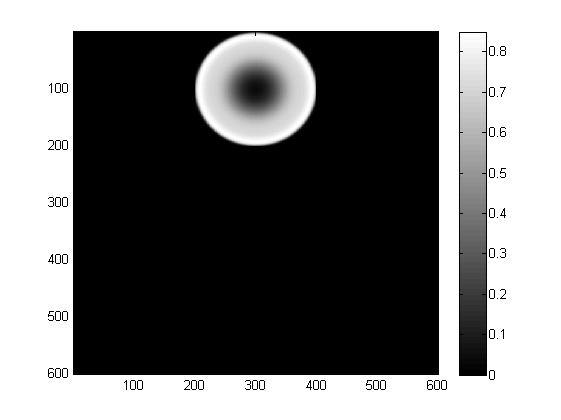}
\includegraphics[scale=0.45]{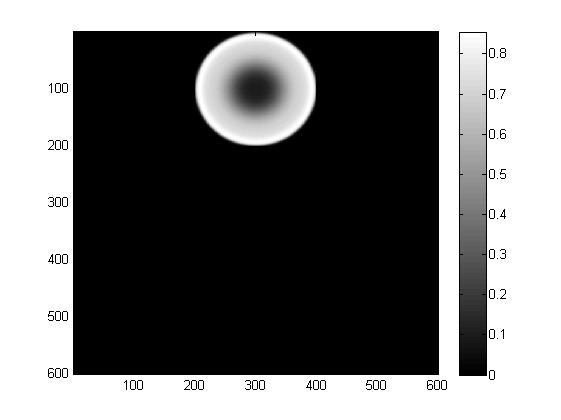}
\caption{We  have presented averaged reconstructions of $f_t$ with $10\%$ and $50\%$ added noise in the left and right hand images respectively.}
\label{figure18}
\end{figure}

We have presented reconstructions of a density which is homogeneous where it is not known to be zero. To give an inhomogeneous example, we have presented reconstructions with varying levels of added noise of a simulated hollow tube cross section in figures \ref{figure17} and \ref{figure18}.

We can summarize our method as follows:
\begin{enumerate}
\item Measure the scattered intensity energy $E_s$ and divide by $P_{\text{avg}}$ and the slice thickness to calculate values for $\mathcal{D}f$.
\item Smooth the data sufficiently and apply approximation (\ref{app1}) to calculate values for $R\tilde{f}$.
\item Reconstruct $\tilde{f}$ by filtered backprojection and recover $f$ from the definition given in Theorem \ref{th1}.
\item Average over a number of source views to improve the image quality and set $f$ to $0$ outside its support.
\end{enumerate}
\clearpage

\clearpage

\section{Conclusion}
\label{conc}

We have proposed a new fast method to determine the electron density in x-ray scanning applications, with a fixed energy sensitive detector machine configuration where it is possible to measure photon intensity in the dark field. We have shown that the density may be reconstructed analytically using the Compton scattered intensity. This method does not require the photon source to be monochromatic as is the case in recent literature, which is important from a practical standpoint as it may not be reasonable to assume a monochromatic source in some applications. Also if the source is monochromatic we cannot gain any insight into the energy dependence of the attenuation coefficient, which would rule out the recent advances in image rendering \cite{pca1, pca2}, where a combination of multivariate and cluster analysis can be used to render a colour x-ray image. 

Using Sobolev space estimates, we have determined an upper bound for the least squares error in our solution in terms of the least squares error in our data. This work is based on the approach taken by Natterer in \cite{nat1}. 

We have shown, under the right assumptions, that the atomic number of the target is determined uniquely by the full data. With this theory in place we intend to pursue a more practical means to reconstruct the atomic number $Z$, as the graph reading method used in the present paper was ineffective in giving an accurate reconstruction for $Z$.

We summarize our method to recover the density image in section \ref{res} and we reconstruct a simulated water bottle cross section via a possible practical implementation of this method. In this simple case the smoothing method (simple moving average) applied was effective and we were able to reconstruct a circular cross section of approximately uniform density. Although in the presence of noise the pixel values of our reconstructed density image on average differed from the original values by as much as $15\%$. We have also provided reconstructions of a simulated hollow tube cross section. In this case the inner edge of the tube cross section appeared quite blurred in the reconstruction when noise was added to the simulated data. We performed a number of trial reconstructions with different randomly generated datasets. The results presented in this paper are typical of our trial results.

We hope also to test our methods through experiment. For example, if we were to take an existing x-ray machine of a similar configuration to that discussed in the present paper, and attach energy sensitive detectors alongside the existing detectors or if we were to replace them, then we could see how closely our forward problem models the intensity of photons measured in the dark field in practice. 

\section*{Acknowledgements}
I would like to thank my Ph.D. supervisor Prof William Lionheart for his guidance and inspiration. The author is also grateful to Prof Robert Cernik for his helpful comments and discussion regarding energy sensitive detectors, and to Dr Ed Morton and Dr Tim Coker of Rapiscan systems for information on baggage scanning. This work has been funded jointly by the EPSRC and Rapiscan systems.

\section*{Appendix -- The RTT80; An example application in threat detection}
The RTT80 (real time tomography) X-ray scanner is a switched source, offset detector CT machine designed with the aim to scan objects in real time. Developed by Rapiscan systems, the RTT80 is currently used in airport security screening of baggage.

\begin{figure}[!h]
\centering
\begin{tikzpicture}[scale=2.9]
\draw [green] (0,0) circle [radius=0.7];
\draw [red, very thin] (-0.3,-0.2) rectangle (0.3,0.2);
\draw [very thin, dashed] (0,-1)--(0,0.7);
\draw [very thin, dashed] (0,-1)--(0.1,0.6928);
\draw [very thin, dashed] (0,-1)--(-0.1,0.6928);
\draw [very thin, dashed] (0,-1)--(0.2,0.6708);
\draw [very thin, dashed] (0,-1)--(-0.2,0.6708);
\draw [very thin, dashed] (0,-1)--(0.3,0.6324);
\draw [very thin, dashed] (0,-1)--(-0.3,0.6324);
\draw [very thin, dashed] (0,-1)--(0.4,0.5744);
\draw [very thin, dashed] (0,-1)--(-0.4,0.5744);
\draw [very thin] (0.9,0.9)--(0.3,0.6324);
\node at (1.1,1) {detector ring};
\draw [blue] (0,0) circle [radius=1];
\draw [very thin] (-1.1,-0.9)--(-0.7,-0.71);
\node at (-1.1,-1) {source ring};
\draw [thick] (-0.4,-0.3)--(0.4,-0.3);
\draw [very thin] (-0.9,0.9)--(-0.3,0);
\node at (-0.9,1) {scanned object};
\draw [very thin] (1.1,-0.9)--(0.3,-0.3);
\node at (1.1,-1) {conveyor belt};
\draw [red, thin] circle [radius=0.57];
\draw [very thin] (1.3,0.3)--(0.57,0);
\node at (1.4,0.4) {scanning tunnel};
\draw [<->, thick] (-1.8,-1)--(-1.8,-0.7);
\draw [<->, thick] (-1.8,-0.7)--(-1.8,-0.57);
\draw [<->, thick] (-1.8,-0.57)--(-1.8,0.57);
\node at (-1.9,-0.85) {1};
\node at (-2,-0.635) {0.375};
\node at (-1.9,0) {5};
\end{tikzpicture}
\caption{The RTT80 machine configuration is displayed. The source-detector ring offset is relatively small and will be modelled as zero. The RTT80's  relative dimensions (the source ring, detector ring and scanning tunnel radii) are presented to the left of the diagram.}
\label{figRTT}
\end{figure}
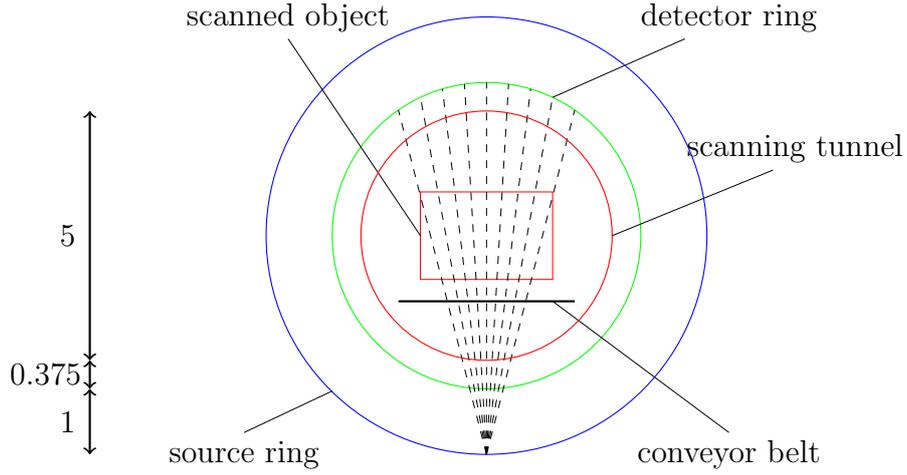

The RTT80 consists of a single fixed ring of polychromatic X-ray sources and multiple offset rings of detectors, with a conveyor belt and scanning tunnel (within which the scanned object would be placed) passing through the centre of both sets of rings. See figure \ref{figRTT}. If the detectors are energy sensitive, then in this case we have the problem of reconstructing a density slice supported within the scanning tunnel from its integrals over toric sections, with tips at the source and detector locations. We wish to check whether it is reasonable to approximate a set of toric section integrals as integrals over discs whose boundaries intersect a given source point, as then we can apply our proposed reconstruction method to reconstruct the density slice analytically. 

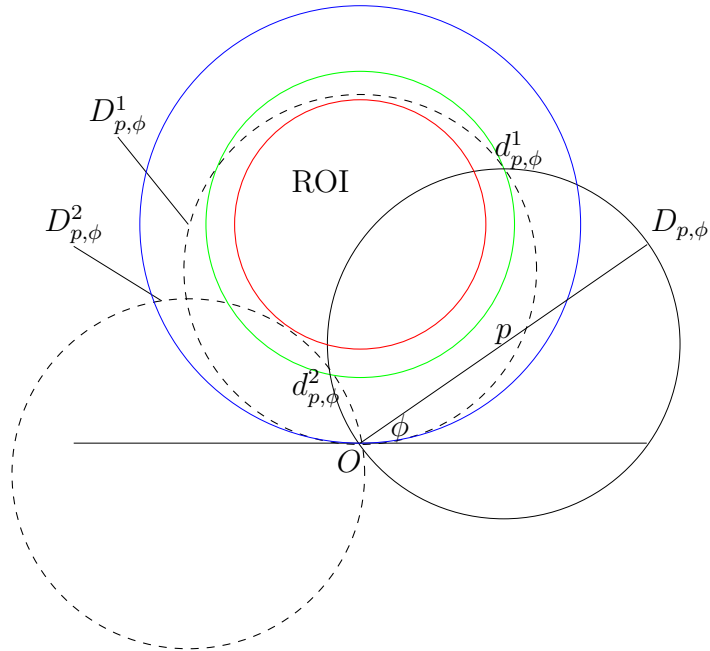
\begin{figure}[!h]
\centering
\begin{tikzpicture}[scale=2.9]
\draw [thin, green] (0,0) circle [radius=0.7];;
\draw [very thin] (0.651,-0.546) circle [radius=0.8];
\node at (-0.05,-1.08) {$O$};
\node at (0.72,0.34) {$d^1_{p,\phi}$};
\node at (0.651+0.8,0) {$D_{p,\phi}$};
\node at (-0.18,0.2) {$\text{ROI}$};
\draw [red, thin] circle [radius=0.57];
\draw [very thin] (-1.3,-1)--(1.3,-1);
\node at (-0.2,-0.75) {$d^2_{p,\phi}$};
\draw [thin, dashed]  (0,-0.206) circle [radius=0.8];
\draw [thin, dashed]  (-0.78,-01.14) circle [radius=0.8];
\draw [very thin] (-1.1,0.4)--(-0.78,0);
\draw [very thin] (-1.3,-0.1)--(-0.9,-0.35);
\node at (-1.1,0.5) {$D^1_{p,\phi}$};
\node at (-1.3,0) {$D^2 _{p,\phi}$};
\draw [very thin,blue] (0,0) circle [radius=1];
\draw [very thin] (0,-1)--(1.302,-0.092);
\node at (0.651,-0.5) {$p$};
\node at (0.18,-0.93) {$\phi$};
\end{tikzpicture}
\caption{The RTT80 configuration is displayed. The origin is denoted by $O$ as in section \ref{sec1}. This is where a source is located. A disc $D_{p,\phi}$ with boundary intersecting $O$ and two detector points $d^1_{p,\phi}$ and $d^2_{p,\phi}$ is shown to have a non empty intersection with the set of points in our ROI (the scanning tunnel). The disc $D^1_{p,\phi}$ is the reflection of $D_{p,\phi}$ in the line segment connecting $O$ to $d^1_{p,\phi}$. Similarly for $D^2_{p,\phi}$.}
\label{figRTT1}
\end{figure}

Let us refer to figure \ref{figRTT1} and let $D_{p,\phi}$ be defined as in section \ref{sec1}. We define the toric sections $T^1 _{p,\phi}=D^1 _{p,\phi}\cap D_{p,\phi}$, $T^2 _{p,\phi}=D^2 _{p,\phi}\cap D _{p,\phi}$, $T^3 _{p,\phi}=D^1 _{p,\phi} \cup D _{p,\phi}$ and $T^4 _{p,\phi}=D^2 _{p,\phi} \cup D_{p,\phi}$. Let $A(S)$ denote the area of a set $S\subseteq \mathbb{R}^2$ and let $T\subseteq \mathbb{R}^2$ denote the set of points within our ROI (region of interest, i.e the scanning tunnel). For a large sample of discs, we will check for every disc $D_{p,\phi}$ in the sample, whether $\exists i\in \{1,2,3,4\}$ such that $A(D_{p,\phi}\cap T)\approx A(T^i _{p,\phi}\cap T)$. 

Let $D_r$ be defined as in section \ref{sec1}. Then if we consider the machine specifications for the RTT80, we can calculate $r=6.75$ and the difference in radius between the detector ring and the scanning tunnel to be $0.375$. See figure \ref{figRTT}. 
For our test, we consider a sample of 36000 discs with diameters $p=1.375+\frac{5(i-1)}{99}$ for $1\leq i\leq 100$ and $\phi=\frac{\pi j}{180}$ for $1\leq j\leq 360$. We have chosen $p\in [1.375,6.375]$ and $\phi \in [0,2\pi]$ values in a range sufficient to determine a unique density slice image for densities supported on $T$. Refer to Corollary \ref{cor1}. For each of our chosen $p$ and $\phi$ value pairs, the difference:
\begin{equation}
\min_{1\leq i\leq 4} \left(A(D_{p,\phi}\cap T)-A(T^i _{p,\phi}\cap T)\right)\approx 10^{-16}
\end{equation}
was found to be negligible. Let $f:\mathbb{R}^2\to \mathbb{R}$ be an example density slice with support contained in $T$. Then for any disc $D_{p,\phi}$ in our sample we have:
\begin{equation}
\label{intD}
\int_{D_{p,\phi}}f=\int_{D_{p,\phi}\cap T}f= \int _{T^i _{p,\phi}\cap T}f=\int_{T^i _{p,\phi}}f
\end{equation}
which holds for some $i\in \{1,2,3,4\}$. So, the integral of $f$ over $D_{p,\phi}$ is equal to at least one of four toric section integrals over $f$. Assuming also that there is little error implied by our physical approximations (these are discussed in detail in section \ref{phys}), the integral (\ref{intD}) would be determined approximately by at least one of four data sets, namely the photon intensity measured for two possible energy levels at two possible detector locations ($d^1_{p,\phi}$ or $d^2_{p,\phi}$). Thus, given that the inverse disc transform is only mildly ill posed (this was determined to be the case in section \ref{sob}, based on the criteria given by Natterer in \cite{nat1}), it seems that we should be able obtain a satisfactory density image reconstruction in this application. 

In airport baggage screening, we are interested in identifying a given material as either a threat or non-threat. Let $n_e$ be the electron density and let $Z$ denote the effective atomic number. We define the threat space to be the set of materials with $(n_e,Z) \in \text{T}$, where $\text{T} \subseteq [0,\infty)\times [1,100]$ is the class of threat $(n_e,Z)$ pairs. For a given suspect material, we can apply the methods presented in this paper to reconstruct $n_e$ and $Z$. Then if $(n_e, Z) \in \text{T}$, we can identify the suspect material as a potential threat. We note that although we failed to obtain an accurate $Z$ reconstruction in the present paper, we aim to show that a more precise determination of $Z$ is possible in future work. Also, the reconstruction methods we have presented should be fast to implement as they are largely based on the filtered back-projection algorithm. This is important in an application such as airport baggage screening, as we require the threat detection method we apply to not only be accurate in threat identification, but to also be an efficient process.  

\clearpage

\appendix


\clearpage

\begin{thebibliography}{10}
\bibitem{pal1}
V.P. Palamodov, ``An analytic reconstruction for the Compton scattering tomography in a plane" Inverse Problems 27 (2011) 125004 (8pp).
\bibitem{cone}
V. Maxim, M. Frandes, R. Prost, ``Analytical inversion of the Compton transform using the full set of available projections" Inverse Problems 25 (2009) 095001 (21pp).
\bibitem{norton}
S.J. Norton, ``Compton scattering tomography" J. Appl. Phys. 76 2007–15 (1994).
\bibitem{pca1}
Q. Xu, H. Yu, J. Bennett, ``Image reconstruction for hybrid true-color micro CT" IEEE Trans Biomed Eng. 59(6) 1711–1719 (2012).
\bibitem{pca2}
C.K. Egan, S.D.M. Jacques, R.J. Cernik, ``Multivariate analysis of hyperspectral hard x-ray images" Wiley 42 151–157 (2012).
\bibitem{he}
Helgason, S. ``Gropes and geometric analysis" Academic Press, Orlando-San Diego-San
Francisco-New York-London-Toronto-Montreal-Sydney-Tokyo-Sao Paulo (1984). 
\bibitem{fit}
Palinkas, G. ``Analytic approximations for the incoherent x-ray intensities of the atoms from Ca to Am" Acta Cryst. A29, 10 (1973).
\bibitem{hub}
Hubbell, J. H. et. al ``Atomic form factors, incoherent scattering functions and photon scattering cross sections" J. Phys. Chem. Ref. Data, Vol. 4, No. 3, (1975).
\bibitem{vieg}
WM. J. Veigele, ``Photon cross sections from 0.1keV to 1MeV for elements $Z=1$ to $Z=94$" Atomic Data Tables, 5, 51-111 (1973).
\bibitem{DJ}
D. F. Jackson, D. J. Hawkes ``An accurate parametrisation of the x-ray attenuation coefficient" Phys. Med. Biol. 25 1167 (1980).
\bibitem{Nic}
Wadeson, N., ``Modelling and correction of scatter in a switched source multi ring detector CT machine" PhD Thesis, University of Manchester UK, (2012).
\bibitem{nat1}
F. Natterer ``The mathematics of computerized tomography" SIAM (2001).
\end{thebibliography}
\end{document}